\documentclass[12pt]{article}
\usepackage{amssymb,amsmath}
\newtheorem{theorem}{Theorem}[section]

\newtheorem{proposition}{Proposition}[section]
\newtheorem{corollary}{Corollary}[section]
\newtheorem{definition}{Definition}[section]
\newtheorem{example}{Example}[section]

\newcommand{\w}{\wedge}
\newcommand{\n}{\notag}
\newcommand{\im}{\mbox{i}}
\newcommand{\hook}{\hookrightarrow}
\newcommand{\tb}{\textbf}
\newcommand{\vsp}{\vspace}
\newcommand{\sq}{$\; \square$}

\newcommand{\E}{ \mathbb{E} }
\newcommand{\C}{ \mathbb{C} }
\newcommand{\R}{ \mathbb{R} }
\newcommand{\Remark}{\emph{Remark.}$\;$}
\newcommand{\Proof} {\emph{Proof.}$\; \,$}

\newcommand{\bomega}{\bar{\omega}}

\newcommand{\bnu}{\bar{\nu}}
\newcommand{\W}{ \mathcal{W} }

\begin{document}

\sloppy

\begin{center}
\Large{\tb{
Willmore Legendrian surfaces in
\\pseudoconformal 5-sphere }}
\end{center}

\begin{center}
Sung Ho Wang\\
\texttt{shw@kias.re.kr}
\end{center}

\vspace{1pc}

0. Introduction

1. Pseudoconformal stereographic projection

2. Legendrian surfaces

3. A duality theorem

4. Willmore Legendrian surfaces with constant dual

Appendix. Euler-Lagrange equation

\vsp{3pc}
\begin{center}
\textbf{\large{Introduction}}
\end{center}

For a compact oriented submanifold $\, X: M^n \hook \E^N$
in a Euclidean space, Willmore functional $\, \W(X)$ is defined by
\begin{equation}
\W(X)=\int_M\, \frac{1}{n} \, |B_0|^n\, dX,\n
\end{equation}
where $\, B_0$ is the trace free part
of the second fundamental form of $\, X$.
It was introduced by Willmore
in a slightly different but equivalent form
for surfaces in $\, \E^3$. The functional is invariant under conformal
transformation, and it is natural to extend it for
submanifolds in conformal $\, N$-sphere $\, S^N=\E^N \cup \{\infty\}$.
Willmore conjecture asks if Clifford torus in $\, S^3$ is
the unique minimizer of $\, \W(X)$ among all immersed tori.

Li and Yau introduced the notion of conformal area of a compact Riemann surface,
and showed that Willmore conjecture is true if the induced complex(conformal)
structure of the competing torus lies in a specific neighborhood of
the square torus [LY].
Bryant proved a duality theorem, and studied a class of Willmore surfaces
for which Willmore dual is constant [Br2][Br3].
He gave a classification of such surfaces in terms of complete minimal surfaces
of finite total curvature with embedded ends of zero logarithmic growth.
Ros answered the conjecture in the affirmative
among tori invariant under antipodal involution using his solution to
isoperimetric problem in $\, \R P^3$ [Ro].
Recently Topping gave an integral geometric proof of Ros's result [To].
Babich and Bobenko constructed examples of Willmore tori in $\, \E^3$ with
a planar reflectional symmetry by gluing two minimal surfaces in hyperbolic
space [BB].
Helein gave a loop group formulation of constructing all Willmore tori [He].
Montiel and Urbano studied a Willmore type functional for
surfaces in $\, \C P^2$ [MU].

The problem has a pseudoconformal(CR) cousin.
Let $\, \Sigma^n =S^{2n+1}$ be the pseudoconformal sphere with
SU$(n+1,1)$ as the group of automorphisms [ChM].
Let $\, X: M^n \hook \Sigma^n$ be a Legendrian submanifold.
Then $\, M$ inherits a conformal structure, and
Willmore functional $\, \W(X)$ in this case is a pseudoconformally invariant
second order functional defined similarly as above, see \eqref{Willf}.
A Legendrian submanifold is called Willmore if it is critical for this functional.

The purpose of this paper is to study Willmore Legendrian surfaces
in the spirit of [Br2].
There are two main results.
The first is a construction of  Willmore Legendrian dual
for a class of Willmore Legendrian surfaces, \tb{Theorem \ref{typeB}}.
The second is a characterization of Willmore Legendrian surfaces
with constant Willmore dual in terms of immersed meromorphic curves in
$\, \C^2$ satisfying an appropriate real period condition,
\tb{Theorem \ref{typeC}}.
We show that \emph{every compact Riemann surface
admits a generally one to one, conformal, Willmore Legendrian
immersion in $\, \Sigma^2=S^5$ with constant Willmore dual},
\tb{Theorem \ref{theorem42}}.
As a corollary,
\emph{
every compact Riemann surface can be conformally immersed in $\, \C^2$
as an exact, algebraic Lagrangian surface.}

  In Section 1, we give an explicit formulation of
pseudoconformal stereographic projection.
The associated pseudoconformal involution
will play a role in our ends analysis of Willmore Legendrian surfaces
with constant Willmore dual.
  In Section 2, a basic structure equation
for Legendrian surfaces in $\, \Sigma^2$  is established.
A complex quartic differential $\, \Phi$ arises as a fourth order
invariant of the immersion.
  In case $\, \Phi=0$, a complex sextic differential $\, \Psi$ is defined
and holomorphic. In addition, the umbilic loci in this case,
the points where the trace free part of second fundamental form vanish,
is a closed subset with no interior.
  In section 3, we define the Willmore functional, and identify its
Euler-Lagrange equation; a Legendrian surface is Willmore
when $\,\Phi$ is holomorphic with respect to the induced complex structure.
It is a fifth order elliptic equation for the immersion.
When $\, \Phi=0$, an associated Willmore dual is well defined and
can be smoothly extended across the umbilic loci.
  In section 4, we give a Weierstra\ss\, type representation
for surfaces with $\, \Phi=\Psi=0$,
which is equivalent to Willmore dual being constant \eqref{data}.
By Riemann-Roch, every Willmore Legendrian sphere belongs to this class.
We show that there exists Willmore sphere with $\, \W = 4\pi k$
for each integer $\, k \geq 0$.

One of our motivation is the question,
"What is the best compact Legendrian surface of given genus
in $\, \Sigma^2=S^5$ ?"
It is perhaps the case that the following is
the pseudoconformal analogue of Willmore conjecture;
\emph{the minimal Legendrian hexagonal torus
\begin{equation}
T=\{ \, (\xi^1, \, \xi^2, \, \xi^3)\, |\; \; |\xi^i|^2= \frac{1}{3}, \;
\mbox{Im} (\xi^1 \xi^2 \xi^3) =0 \, \} \subset S^5 \subset \C^3\n
\end{equation}
uniquely minimizes Willmore functional among all immersed Legendrian tori.}
A study of this and other related problems for Willmore Legendrian tori
will be the subject of the subsequent paper.

\vsp{2pc}
\section{Pseudoconformal stereographic projection}
In this section, we describe 5-sphere $\, S^5$ as the standard
model of pseudoconformal space. Pseudoconformal analogue of
the usual stereographic projection in conformal geometry is introduced, which
enables us to compare the geometry of Legendrian surfaces in $\, S^5$
with the geometry of exact Lagrangian surfaces in $\, \mathbb{C}^2$.

Let $\, \mathbb{C}^{3,1}$ be the complex vector space with coordinates
$\, z = (\, z^0, \, z^i, \, z^{3} \,)$, $\, 1 \, \leq i \, \leq  2$,
and a Hermitian scalar product
\begin{equation}
\langle \, z, \, \bar{z} \, \rangle = z^i \, \bar{z}^i +
\mbox{i} \, (z^0 \, \bar{z}^3 - z^3 \, \bar{z}^0 ).\n
\end{equation}
Let $\, \Sigma$ be the set of equivalence classes up to scale of
null vectors with respect to this product. Let SU$(3,1)$ be the
group of unimodular linear transformations that leave the form
$\, \langle \, z, \, \bar{z} \, \rangle$ invariant. SU$(3,1)$ acts
transitively on $\, \Sigma$, and
\begin{equation}\label{proj}
p: \, \mbox{SU}(3,1) \to \Sigma = \, \mbox{SU}(3,1)/P
\end{equation}
for an appropriate  subgroup $\, P$ [ChM].

Explicitly, consider an element $\, Z=(\, Z_0, \, Z_i, \, Z_{3}\,)
\in \mbox{SU}(3,1)$ as an ordered set of  four column vectors
in $\, \mathbb{C}^{3,1}$  such that det$(Z)=1$, and that
\begin{equation}\label{prod}
\langle \, Z_i, \, \bar{Z}_j  \, \rangle = \, \delta_{ij}, \quad \langle \,Z_0,
\, \bar{Z}_{3} \, \rangle = - \langle \,Z_{3}, \, \bar{Z}_0 \, \rangle=\, \mbox{i},
\end{equation}
while all other scalar products are zero. We define
\begin{equation}
p(Z) \, = \, [Z_0],  \n
\end{equation}
where $\, [Z_0]$ is the equivalence class of null vectors
represented by   $\, Z_0$. The left invariant Maurer-Cartan form $\, \pi$
of SU$(3,1)$ is defined by the equation
\begin{equation}
dZ = Z \, \pi,\n
\end{equation}
which is in coordinates
\begin{equation}\label{struct1}
d (\, Z_0, \, Z_i, \, Z_{3} \, ) = ( \, Z_0, \, Z_j, \, Z_{3} \, ) \,
\begin{pmatrix}
\pi_0^0 & \pi_i^0 & \pi_{3}^0 \\
\pi_0^j & \pi_i^j & \pi_{3}^j \\
\pi_0^{3} & \pi_i^{3} & \pi_{3}^{3}
\end{pmatrix}.
\end{equation}
Components of $\, \pi$ are subject to the relations obtained from
differentiating (\ref{prod}) which are
\begin{align}
\pi_0^{3}&= \bar{\pi}_0^{3}, \quad \pi_{3}^{0}= \bar{\pi}_{3}^{0} \label{strela} \\
\pi_i^{3}&= - \, \mbox{i} \, \bar{\pi}_0^i, \quad \pi_{3}^i = \mbox{i} \, \bar{\pi}_i^0 \notag \\
\pi_j^i + \bar{\pi}_i^j &=0, \quad \pi_0^0 + \bar{\pi}_{3}^{3} =0  \notag \\
\mbox{tr} \,  \pi &=0,  \n
\end{align}
and $\, \pi$ satisfies the structure equation
\begin{equation}
- \, d\pi = \pi \wedge \pi.\n
\end{equation}

It is well known that the SU$(3,1)$-invariant CR structure on
$\, \Sigma \subset \mathbb{C}P^{3}$ as a real hypersurface is
biholomorphically equivalent to the standard  CR structure on
$\, S^{5} = \partial \tb{B}^{3}$, where
$\, \tb{B}^{3} \subset \mathbb{C}^{3}$ is the unit ball.
The structure equation (\ref{struct1}) shows that for any local section
$\, s: \, \Sigma \, \to$ SU$(3,1)$, this CR structure is defined by
the hyperplane fields $\, (s^*\pi_0^{3})^{\perp} = \mathcal{H}$
and the set of (1,0)-forms $\, \{ \, s^*\pi_0^i \, \}$.

Let $\, \mathbb{L} \subset \mathbb{C}^{3,1}$
be the cone of nonzero null vectors.
SU$(3,1)$ acts transitively on $\, \mathbb{L}$, and
$\, \mathbb{L} \to \Sigma$ is a $\, \mathbb{C}^*$-bundle.
Take $\, w_3 \in \mathbb{L}$, and let
\begin{equation}
\Sigma_{w_3} = \, \{ \, w \in \mathbb{L} \; | \; \; \langle \, w, \, \bar{w}_3
                      \, \rangle = \im \, \}. \n
\end{equation}

As a submanifold of $\, \mathbb{C}^{3,1}$, $\, \Sigma_{w_3}$ inherits
a real contact form $\, \im \, \langle \, dw, \, \bar{w} \, \rangle$ and
a degenerate Hermitian metric of signature (2,0), which together  form
a subHermitian structure. Moreover, the natural projection
$\, \pi: \, \Sigma_{w_3} \to \Sigma - \{ p_{\infty} \}$ is easily seen to
be a pseudoconformal equivalence, where $\, p_{\infty} =[w_3]$.

Take $\, w_0 \in \Sigma_{w_3}$, and denote
$\, \mathbb{E} = \, \langle w_0, \, w_3 \rangle^{\perp}
\simeq \mathbb{C}^2 \subset \mathbb{C}^{3,1}$.
From the choice of $\, w_0$, the induced Hermitian metric on $\, \mathbb{E}$
is positive definite. Let $z = (z^1, \, z^2)$ be the standard coordinate
of $\, \mathbb{E}$.

Let $\, P: \Sigma_{w_3} \to \mathbb{E}$ be the orthogonal projection
defined by
\begin{equation}\label{sprojection}
P(w)= w - w_0 - \im \, \langle \, w, \, \bar{w}_0 \, \rangle \, w_3.
\end{equation}
Let $\, \dot{w}$ be a tangent vector to $\, \Sigma_{w_3}$ so that
$\, \langle \, \dot{w}, \, \bar{w}_3 \, \rangle =0$. Then
\begin{equation}
\langle \, P_*(\dot{w}), \, \bar{P}_*(\dot{w}) \, \rangle
 = \langle \, \dot{w}, \, \bar{\dot{w}} \, \rangle,  \n
\end{equation}
and $\, P$ is  a Riemannian submersion whose fibers are in the null
direction. In fact, since $\, P$ is linear, it is easy to see that
\begin{equation}
P(w) = P(w') \iff w - w' = s \, w_3, \quad s \in \mathbb{R},  \n
\end{equation}
and $\, \Sigma_{w_3}$ is ruled by real lines parallel to $\, w_3$.

\begin{definition}
Let $\, p_{\infty} \in \Sigma$, and let $\, w_0, \, w_3 \in \mathbb{L}$
be such that $\, p_{\infty} =[w_3]$,
and that $\, \langle  w_0,  \, \bar{w}_3 \rangle = \textnormal{\im}$.
The pseudoconformal stereographic projection at $\,  p_{\infty}$ is
the (orthogonal) projection defined by \eqref{sprojection}
\begin{equation}
P: \,  \Sigma - \{p_{\infty}\} \simeq \Sigma_{w_3} \to \mathbb{E}
   = \, \langle w_0, \, w_3 \rangle^{\perp} \simeq \mathbb{C}^2.  \n
\end{equation}
\end{definition}

\emph{Remark.}$\;$ Given a point $\,  p_{\infty} \in \Sigma$, the stereographic
projection from $\,   p_{\infty}$ depends on the choice
null vectors $\, w_0, \, w_3 \in \mathbb{L}$.
For definiteness, we assume such a choice is made whenever stereographic
projection is used.

Let $\,t$ be the coordinate of $\, \mathbb{R}$, and
let $\, V = \mathbb{R} \oplus \mathbb{E}$ with the \emph{standard}
subHermitian structure given by the contact form
$\, dt + \frac{\im}{2} \, ( \, \langle \, dz, \, \bar{z} \, \rangle
                            -  \langle \, z, \, d\bar{z} \, \rangle \, )$,
and $\, \E$-wise subHermitian metric.

\begin{proposition}
There exists a quadratic subHermitian isomorphism
\begin{equation}
\hat{Q}: \, V \, \to \,  \Sigma_{w_3}. \n
\end{equation}
\end{proposition}
\emph{Proof.}$\;$  Let $\, Q: \mathbb{E} \to  \Sigma_{w_3}$ be the quadratic
map
\begin{equation}
Q(z) = z+w_0 - \frac{\im}{2} \langle \, z, \, \bar{z} \, \rangle \, w_3 \n
\end{equation}
for $\, z \in \mathbb{E}$. Then
$\, \langle Q(z), \, \bar{w}_3 \rangle = \im$,
$\, \langle Q(z), \, \bar{Q}(z) \rangle = 0$, and $\,Q$ is well defined.

Note
\begin{align}\label{Qinverse}
P(Q(z)) &= z \\
\langle Q_*(\dot{z}), \, \bar{Q}_*(\dot{z}) \rangle &=
\langle  \dot{z}, \, \bar{\dot{z}} \rangle, \n
\end{align}
for $\, \dot{z} \in \mathbb{E}$. Now define
\begin{equation}\label{Qform}
\hat{Q}(t, \, z) = Q(z) + t\, w_3.
\end{equation}
From \eqref{Qinverse}, $\, Q$ is one to one. Since fibers of $\, P$ are real
lines parallel to $\, w_3$, $\, \hat{Q}$ is also one to one, and obviously onto,
and hence a diffeomorphism. The rest follows by direct computation.
The inverse of $\, \hat{Q}$ is given by
\begin{equation}
\hat{P}(w) =\hat{Q} ^{-1}(w) =
\big(\frac{\im}{2} ( \langle w, \, \bar{w}_0 \rangle -
\langle w_0, \, \bar{w} \rangle), \; P(w) \big). \quad \square \n
\end{equation}

$\Sigma = V \cup \{p_{\infty}\}$ can thus be identified with the one point
pseudoconformal compactification of $\, V$, and Legendrian
surfaces in $\Sigma - \{p_{\infty}\}$ are in one to one correspondence
with \emph{exact} Lagrangian surfaces in $\, \mathbb{C}^2$ via
pseudoconformal stereographic projection. A Lagrangian surface
$\, M \subset \, \mathbb{C}^2$ with symplectic form $\, \varpi$ is exact
if an(any) antiderivative $\, \eta$ of $\, \varpi$ in $\, \C^2$,
$\, d \eta = \varpi$, is an exact 1-form when restricted to $\, M$.

For later purpose, we describe the pseudoconformal involution.
Suppose we swap $\, (w_3, \, w_0)$ to $\, (w_0, \, -w_3)$ in the above
construction. Under the isomorphism \eqref{Qform}, this induces
the pseudoconformal involution of $\, \Sigma = V \cup \{ p_{\infty} \}$
given by
\begin{equation}\label{involution}
(\, t, \, z) \to (-\mbox{Re}\, \lambda, \, \lambda \, z),
\quad \lambda = \frac{1}{t-\frac{\im}{2}\langle z, \,\bar{z} \rangle}.
\end{equation}
It interchanges $\, p_{\infty}=[w_3]$ with
$\, (0, \,0) = [w_0] \in V$.

Let $\, M$ be a compact surface with a divisor
$\, D \subset M$. Let $\, f: M-D \hook \mathbb{E}$ be an exact Lagrangian
surface which is complete and with embedded ends $\, D$.
If $\, f$ satisfies an appropriate asymptotic regularity condition
on $\, D$, the natural Legendrian lift $\, X_f$ of $\, f$ in
$\, V  = \R \oplus \E \simeq \Sigma_{w_3}$ closes up across $\, D$
to be a smooth Legendrian immersion $\, X_f: M \hook  \Sigma$.
One such occasion will be described in \tb{Section 4}.

\vsp{2pc}
\section{Legendrian surfaces}
To study Legendrian surfaces in $\, \Sigma$,
we employ the method of moving frames and construct
an adapted frame bundle, on which the rest of our analysis is based.
Willmore functional for Legendrian surfaces we are interested in
is nothing but the well known Willmore functional
restricted to exact Lagrangian surfaces in $\, \mathbb{C}^2$
under pseudoconformal stereographic projection.

The projection map \eqref{proj} makes SU$(3,1)$ into
a right principal bundle over
pseudoconformal 5-sphere $\, \Sigma$ with fiber
\begin{equation}
G_0 = \left\{
\begin{pmatrix}
a^{-1}  &   \im \, a^{-1} \, \bar{v}^t \, A   &  c  \\
0       &    A  &  v               \\
0       &    0  &  \bar{a}
\end{pmatrix} \Bigg| \quad
\begin{array}{ll}
&a\ne 0, \, v^t = (v^1, \, v^2),
\, |v|^2=\im (\bar{a}\bar{c}-ac)\\
&\mbox{$\, \, A \in U(2), \,a^{-1} \bar{a}$ det$\, A = 1$}
\end{array}
\right\}.\n
\end{equation}
This verifies that the 1-form $\, \pi_0^3$ from \eqref{struct1} is well defined
up to scale on $\, \Sigma$, and defines a contact structure
$\, \mathcal{H} = (\pi_0^{3})^{\perp}$.
The 1-forms $\, \{ \, \pi_0^1, \, \pi_0^2 \, \}$ are well defined
mod $\, \pi_0^3$ up to conformal unitary transformations,
and defines a conformal Hermitian structure on $\, \mathcal{H}$.

A \emph{Legendrian surface} consists of an oriented surface $\, M$
and a Legendrian immersion
\begin{equation}
X: \, M \hook \Sigma, \quad X_*(TM) \subset \mathcal{H}. \n
\end{equation}
We wish to study its properties invariant under the group action
\eqref{proj}.

The 0-adapted frame bundle $\, F_0$ of $\, X$ is by definition
\begin{equation}\label{F0}
F_0 = X^*\, \mbox{SU(3,1)}.
\end{equation}
Since $\, X$ is Legendrian,
$\, \pi_0^3 = 0$, $\, \pi^i_0 \w \bar{\pi}^i_0=0$ on $\, F_0$, and
$\, \{ \, \pi_0^1, \, \pi_0^2 \, \}$ are semi-basic 1-forms on $\, M$.
An element $\, g \in G_0$ acts on $\, \{ \, \pi_0^1, \, \pi_0^2 \, \}$ by
\begin{equation}
R_g^*
\begin{pmatrix}
\pi_0^1 \\
\pi_0^2
\end{pmatrix} = a^{-1} \, A^{-1} \,
\begin{pmatrix}
\pi_0^1 \\
\pi_0^2
\end{pmatrix},\n
\end{equation}
where $\, R_g$ is the right action by $\, g$.

We may thus define the 1-adapted frame bundle $\, F_1$ by
\begin{equation}
F_1 = \{ \,  \pi_0^i = \omega^i \; \mbox{is real}, \,
                           \omega^1 \wedge \omega^2 > 0 \, \}\subset F_0,\n
\end{equation}
with the structure group
$\,
G_1 = \{ \, g \in G_0 \, | \,a  A = \bar{a} \bar{ A}, \, a^4>0 \, \}.
$
Note that we have used the fact $\, M$ is oriented.
An oriented Legendrian surface inherits a conformal(complex) structure,
for $\, (aA)\, (aA)^t = a\, \bar{a}$, det$\, (aA)= \frac{a^4}{a\,\bar{a}}>0$.

The structure equation \eqref{struct1} with this relation gives
\begin{equation}
 - d \omega^i = ( \pi^i_j - \, \delta_{ij} \, \pi^0_0 ) \wedge \omega^j. \n
\end{equation}
Set $\,  \pi^i_j - \, \delta_{ij} \, \pi^0_0 = \alpha^i_j+\im \, \beta^i_j$.
Then
\begin{align}
- d \omega^i &= \alpha^i_j \wedge \omega^j,  \n\\
     0       &= \beta^i_j \wedge \omega^j.\n
\end{align}

From the relations \eqref{strela}, $\, \pi^i_j+\bar{\pi}^j_i=0$ implies
\begin{align}
\alpha^i_j + \alpha^j_i &= - \delta_{ij} \, (\pi^0_0 + \bar{\pi}^0_0)
         = - \delta_{ij}\,  2 \, \phi, \quad \phi \; \mbox{ is real}, \n \\
\beta^i_j - \beta^j_i   &= 0,\n
\end{align}
so that in particular
\begin{align}
\mbox{tr}\, (\alpha^i_j) &= -2 \, \phi, \n\\
\im \, \mbox{tr}\, (\beta^i_j) &= -2 \, ( \pi^0_0 - \bar{\pi}^0_0).\n
\end{align}
By Cartan's lemma, $\, \beta^i_j = h_{ijk} \, \omega^k$ for a
coefficient $\, h_{ijk}$ that is fully symmetric in its indices.
A computation shows that an element $\, g \in G_1$ now acts by
\begin{equation}
R_g^* \, \im \, ( \pi^0_0 - \bar{\pi}^0_0)  =
\im \, ( \pi^0_0 - \bar{\pi}^0_0)+
(a^{-1}\bar{v}^1+\bar{a}^{-1}v^1) \, \omega^1+
(a^{-1}\bar{v}^2+\bar{a}^{-1}v^2) \, \omega^2.\n
\end{equation}

We may thus define the 2-adapted frame bundle $\, F_2$ by
\begin{equation}
F_2 = \{ \, \mbox{tr}\, (\beta^i_j) = 0 \, \}\subset F_1,\n
\end{equation}
with the structure group
$
G_2 = \{ \, g \in G_1 \, | \,  a\,v+\bar{a}\, \bar{v}=0  \, \}.
$

The geometric meaning of this normalization is clear.
Let $\, p_{\infty}(m)= [Z_3(m)] \in \Sigma$ for $\, m \in M$,
and let $\, P = \Sigma - \{p_{\infty}(m)\} \to \mathbb{C}^2$ be
the stereographic projection at $\, p_{\infty}(m)$.  For any 2-adapted
frame $\, Z = (Z_0, \, Z_i, \, Z_3)$, the Lagrangian surface
$\, P \circ X: M - X^{-1}(p_{\infty}(m)) \hook \mathbb{C}^2$ has zero
mean curvature vector at $\, m$.

Set $\, \pi_3^i = \mu^i + \im \, \nu^i$. Differentiating
$\,  \mbox{tr}\, (\beta^i_j) = 0$, we get $\, \mu^k \wedge \omega^k = 0$.
By Cartan's lemma, $\, \mu^i = p_{ij} \, \omega^j$ for a coefficient
$\, p_{ij} = p_{ji}$. Using the group action by $\, G_2$ as before,
we may translate $\,  \mbox{tr}\, (p_{ij}) =0$.
We may thus define the 3-adapted frame bundle $\, F$ by
\begin{equation}\label{F3}
F = \{ \,   \mbox{tr}\, (p_{ij}) = 0 \, \}\subset F_2.
\end{equation}
Differentiating $\, \mbox{tr}\, (p_{ij}) = 0$, we get
$\, \pi_3^0 = 2\,q_i \, \omega^i$ for a coefficient $\, q_i$.

Since we'll utilize complex geometry later on, we introduce
complex notation here. Set $\, Z_{\pm} = Z_1 \pm \, \im \, Z_2$, and denote
\begin{align}
\omega  &= \omega^1 + \, \im \,\omega^2, \; \rho = \alpha^1_2, \;
\mu     = \mu^1    - \, \im \,\mu^2, \;
\nu     = \nu^1    + \, \im \,\nu^2, \n \\
h &= h_{111} - \im \,h_{112}, \; p = p_{11} - \im \,p_{12}, \;
q = q_{1} - \im \,q_{2},\n
\end{align}
so that
\begin{equation}
\beta^1_1 - \im \beta^1_2 = h \, \omega, \;
\mu = p \, \omega, \; \pi_3^0 = q \, \omega + \bar{q} \, \bar{\omega}.\n
\end{equation}
We summarize the results of the analysis as follows.
\begin{proposition}
Let $\, X: \,M \hook \Sigma$ be an oriented Legendrian surface.
Then there exists a 3-adapted frame bundle
$\, F \subset X^*$SU(3,1) on which
the following structure equations hold.
\begin{equation}\label{struct9}
d(Z_0, \, Z_+, \, Z_-, \, Z_3) =  (Z_0, \, Z_+, \, Z_-, \, Z_3) \,
\begin{pmatrix}
\phi & \nu + \im \, \bar{p} \, \bomega
& \bnu + \im \, p \, \omega & q\, \omega + \bar{q} \, \bomega \\
\frac{1}{2} \, \bomega & \im \, \rho & \im \, h \, \omega
& \frac{1}{2}( \im \, \bnu +  p \, \omega) \\
\frac{1}{2} \, \omega & \im \, \bar{h} \, \bomega & - \im \, \rho
& \frac{1}{2}( \im \, \nu +  \bar{p} \, \bomega) \\
0 & - \im \, \omega & - \im \, \bomega & - \phi
\end{pmatrix}.
\end{equation}
The structure coefficients $\, h, \, p, \, q$ satisfy
the following equations.
\begin{align}
dh + (3 \im  \, \rho +\, \phi)\, h &=\quad \quad \quad \quad \quad\,
- p \,\bomega + z \, \omega, \label{struct90}\\
dp + (2 \im \, \,\rho + 2 \, \phi)\, p&=
- h \, \nu  \quad \quad \; \; \; - q\, \bomega + y \, \omega, \n\\
dq + (\im \, \rho +3 \, \phi)\, q &=
- p\, \nu  \quad \quad \; \; \; - r \,\bomega + x\, \omega, \n\\
dz + (4 \im  \, \rho + 2  \, \phi)\, z&\equiv
- 2 h \, \bnu  \quad \quad \; - z_{-1}\,\bomega \quad \quad \mod \omega, \n\\
dy +(3\im \, \rho+3\, \phi)\, y &\equiv
-2p\, \bnu -z\, \nu -y_{-1}\,\bomega \quad \quad \mod \omega, \n\\
dx +(2\im \, \rho+4\, \phi)\, x&\equiv
-2q\, \bnu - y\, \nu -x_{-1}\,\bomega \quad \quad \mod \omega, \n
\end{align}
where $\, z, \, y, \, x, \, x_{-1}$ are complex coefficients,
$\, r$ is real, and
\begin{equation}
z_{-1} = y+3 h\, |h|^2 , \; y_{-1} = x + 3 |h|^2  \, p.\n
\end{equation}
\end{proposition}
\Remark \, When $\, r=0$, $\,x_{-1} = |h|^2  \, q + p^2 \, \bar{h}$.

Set
\begin{equation}\label{phi}
\Phi = (hq-\frac{1}{2} \, p^2)\, \omega^4.
\end{equation}
A computation with the structure equation shows that  $\, \Phi$ is
a well defined complex quartic differential on $\, M$.
Suppose $\, \Phi \equiv 0$ on $\, M$, and consider
\begin{equation}\label{psi}
\Psi = (zx-\frac{1}{2} \, y^2)\, \omega^6.
\end{equation}
The structure equation shows that
$\, r \equiv 0, \, hz-py+qz \equiv 0$,
and that $\, \Psi$ is a well defined complex
sextic differential on $\, M$ that is holomorphic under
the induced complex(conformal) structure.
This observation will be important in our characterization
of a class of Legendrian surfaces in \tb{Section 4}.

We close this section with an observation on the \emph{umbilic locus}
$\, \mathcal{U}_X = \{ \, m \in M\, | \; h(m) =0 \, \}$.
\begin{proposition}\label{umbilic}
Let $\, X: M \hook \Sigma$ be a connected Legendrian surface with
$\, \Phi \equiv0$.
Then either $\, h \equiv 0$, or the umbilic locus $\, \mathcal{U}_X$ is
a closed subset with no interior.
\end{proposition}
\Proof
Let $\, m \in \mathcal{U}_X$, and let $\, \xi$ be a local coordinate
of $\, M$ centered at $\, m$. Since every conformal structure on a surface
is locally trivial, take a section of the 3-adapted bundle $\, F$ such that
$\, \rho= \, \phi = 0$ [Br2, p41].

First note that $\,  \partial_{\bar{\xi}} \, (hq- \frac{1}{2}p^2)\equiv 0$
implies $\, h r \equiv 0$. Suppose
$\, r(m) \ne 0$. Then $\, h =0$ on a neighborhood of $\, m$, which implies
$\, r=0$ by \eqref{struct90}, a contradiction. Hence $\, r \equiv 0$, and
from \eqref{struct90} again, we compute
\begin{equation}
\partial_{\bar{\xi}}
\begin{pmatrix}
h \\ p \\ q
\end{pmatrix} =
\begin{pmatrix}
0 & -1 & 0 \\
-\nu(\partial_{\bar{\xi}}) & 0 & -1 \\
0 & -\nu(\partial_{\bar{\xi}}) & 0
\end{pmatrix}
\begin{pmatrix}
h \\ p \\ q
\end{pmatrix}.\n
\end{equation}
By a well known theorem [Br2, p42], there exists an integer $\, k(m)$,
$\, 0 \leq k(m) \leq \infty$, such that
\begin{equation}
\begin{pmatrix}
h \\ p \\ q
\end{pmatrix} = \xi^k
\begin{pmatrix}
h_1 \\ p_1 \\ q_1
\end{pmatrix},\n
\end{equation}
for some smooth coefficients $\, h_1, \, p_1, \, q_1$ such that
$\, (h_1(m), \, p_1(m), \, q_1(m)) \ne (0, \, 0, \, 0)$.
Assume $\, k < \infty$ first. Note
$\, h_1q_1 - \frac{1}{2}p_1^2 =0$, and these coefficients also satisfy
\begin{equation}
\partial_{\bar{\xi}}
\begin{pmatrix}
h_1 \\ p_1 \\ q_1
\end{pmatrix} =
\begin{pmatrix}
0 & -1 & 0 \\
-\nu(\partial_{\bar{\xi}}) & 0 & -1 \\
0 & -\nu(\partial_{\bar{\xi}}) & 0
\end{pmatrix}
\begin{pmatrix}
h_1 \\ p_1 \\ q_1
\end{pmatrix}.\n
\end{equation}

Suppose $\, h_1(m) \ne 0$, then $\,m$ is obviously an isolated zero of $\, h$.

Suppose $\, h_1(m) =0, \,p_1(m) \ne 0$. Then
$\, \partial_{\bar{\xi}} \,h_1(m) \ne 0$, and
$\, m$ is not an interior point of the zero set of $\, h_1$,
hence of $\, h = \xi^k \, h_1$.

Suppose $\, h_1(m) =0,\, p_1(m) = 0, \,  q_1(m) \ne 0$.
Then as above, $\,\partial_{\bar{\xi}} \,p_1(m) = -q_1(m)\ne 0$,
and $\, m$ is not an interior point of the zero set of $\, p_1$.
Since $\, h_1 = \frac{1}{2q_1} p_1^2$ in a neighborhood of $\, m$,
$\, m$ is not an interior point of the zero set of $\, h_1$,
hence of $\, h = \xi^k \, h_1$.

The argument above shows that each of the sets
$\,
M_0 = \{ \, m \, | \,k(m) < \infty \}$,
$M_{\infty} = \{ \, m \, | \,k(m) = \infty \}
$,\,
is open. The proposition follows for $\, M$ is connected.
\sq

\vsp{2pc}
\section{A duality theorem}

We continue to use the notations in \tb{Section 1} and \tb{Section 2}.

Let
\begin{equation}\label{Willf}
\Omega_X = \im \, |h|^2 \,  \omega \w \bomega.
\end{equation}
The structure equations \eqref{struct9}, \eqref{struct90} shows
this 2-form is well defined on $\, M$.
\begin{definition}
Willmore functional for a compact   Legendrian surface
$\, X: \,M \hook \Sigma$ is the second order functional
\begin{equation}
\mathcal{W}(X) = \int_M \, \Omega_X.\n
\end{equation}
\end{definition}

Let $\, p_{\infty} \in \Sigma - X(M)$, and
$\, P: \Sigma - \{ p_{\infty}\} \to \C^2$
be the stereographic projection from $\, p_{\infty}$.
Then it is easy to check
\begin{equation}\label{functional}
\W(X) = \int_M \, \frac{1}{2} |B_0|^2 dX,
\end{equation}
where $\, B_0$ is the trace free part of the second fundamental form
of the Lagrangian surface $\, P \circ X(M)$.

A problem of interest of course is to find, if exists, the oriented
Legendrian surface that minimizes the functional among the set of
compact Legendrian surfaces of fixed genus $\, g$.
For example in genus 0 case,
the absolute minimum value 0 is attained
by totally geodesic Legendrian 2-spheres,
i.e., Legendrian surfaces for which $\, h \equiv 0$.

\Remark Minicozzi in his thesis studied the problem of minimizing a functional
equivalent to \eqref{functional} among Lagrangian tori in $\, \C^2$,
and proved the existence and regularity of the embedded minimizer [Min].
Note however an embedded Lagrangian surface in $\, \C^2$ cannot be exact
due to the existence of holomorphic disk spanning the Lagrangian surface
[Gro], and the present geometric variational problem is, although related,
essentially distinct from Minicozzi's work.

A Legendrian surface $\, X: \,M \hook \Sigma$ is called
\emph{Willmore} if it is critical for $\, \W$ for any compactly
supported Legendrian variation. The purpose of this section is to identify
the Euler-Lagrange equation for this variational problem,
and analyze its geometric consequences.

\begin{proposition}\label{EL}
A Legendrian surface $\, X: \,M \hook \Sigma$ is Willmore if the complex
quartic form $\, \Phi$, \eqref{phi}, is holomorphic.
\end{proposition}
We postpone the proof to the appendix.
From the structure equation \eqref{struct90}, this is equivalent to
$\, r = 0, $
which we assume from now on.

\begin{example}\textnormal{
Let $\, X: M \hook S^5=U_3/U_2$ be an oriented Legendrian surface.
From the general theory moving frames,
the second fundamental form of $\, X$
can be expressed as a symmetric cubic differential
$\, \beta \in C^{\infty}(S^3 \,T^*M)$.
Let $\, \beta_0$ be the trace free part of $\, \beta$,
and let $\, 2 \eta = \frac{1}{2}$tr$\, \beta$
be the mean curvature 1-form of $\, X$.
For a local oriented orthonormal coframe $\, \{ \omega^1, \, \omega^2 \}$,
we have
\begin{align}
\beta_0 &= \mbox{Re} \, (h_1-\mbox{i}\, h_2)
  (\omega^1+ \mbox{i}\,\omega^2)^3 \n \\
\eta &= \delta_1 \omega^1 + \delta_2 \omega^2, \n
\end{align}
for coefficients $\, h_i, \, \delta_i$.
The Gau\ss  \, curvature of the induced metric of $\, X$
for instance is $\, K = 1 - 2(h_1^2+h_2^2)+2(\delta_1^2+\delta_2^2)$.
The Willmore functional is written as
\begin{equation}
\mathcal{W}(X) = \int_M \, 2(h_1^2+h_2^2) \, \omega^1 \w \omega^2.\n
\end{equation}
A computation shows that $\, X$ is Willmore when
\begin{equation}\label{willmores5}
*d* ( \frac{1}{2} \Delta \eta + \eta) +
2 \nabla_{\eta^{\sharp}} |\eta|^2 +
\frac{2}{3} \langle \eta^{\sharp} \lrcorner \, \beta_0,
\, \nabla \eta \rangle =0.
\end{equation}
Here * is the Hodge dual operator,
$\, \eta^{\sharp} \in C^{\infty}(TM)$ is
the metric dual of the 1-form $\, \eta$,
and $\, \nabla \eta \in C^{\infty}(S^2T^*M)$
is the covariant derivative of $\, \eta$ (note $d\eta=0$).
Let $\, D \delta_i = \delta_{ij} \, \omega^j$ be the covariant
derivatives of $\, \delta_i$.
Then
\begin{equation}
\frac{1}{3} \langle \eta^{\sharp} \lrcorner \, \beta_0,
\, \nabla \eta \rangle =
( \delta_{11} - \delta_{22})( \delta_1 h_1 + \delta_2 h_2)
+2 \delta_{12} ( \delta_1 h_2 - \delta_2 h_1). \n
\end{equation}
When $\, \eta=0$ and $\, X$ is minimal,
the quartic differential \eqref{phi} vanishes.}
\end{example}

From the arguments in \eqref{phi}, \eqref{psi},
Willmore Legendrian surfaces can be naturally divided into
three disjoint classes.

\qquad \qquad \qquad \tb{A}. $\, \Phi \not \equiv 0$:
general case

\qquad \qquad \qquad \tb{B}. $\, \Phi \equiv 0, \, \Psi \not \equiv 0$:
Willmore dual(to be explained below) is defined

\qquad \qquad \qquad \tb{C}. $\, \Phi \equiv 0, \, \Psi \equiv 0$:
Willmore dual, when it is defined, is constant

Note when $\, M = S^2$, every Willmore Legendrian immersion is of type
\tb{C} due to Riemann-Roch theorem.
The type \tb{C} surfaces will be completely
characterized in terms of Weierstra\ss  \, type holomorphic data in
\tb{Section 4}.

The main result of this section is the following
duality theorem for type \tb{B} surfaces analogous to Bryant's
duality theorem [Br2, Theorem C].
\begin{theorem}\label{typeB}
Let $\, X: \,M \hook \Sigma$ be a Willmore Legendrian surface of type
\tb{B}. Then there exists a smooth Willmore dual $\, \hat{X}: \,M \hook \Sigma$,
which is a weakly conformal, generally singular, Willmore Legendrian immersion
of type \tb{B}. The singular locus of $\, \hat{X}$, the points where
$\, \hat{X}$ is not an immersion, is a closed set with no interior.
$\, \hat{\hat{X}}=X$ whenever it is defined.
\end{theorem}
\Proof  From the structure equations \eqref{struct9}, \eqref{struct90},
consider the map $\, Y: F \to \mathbb{L}$ defined by
\begin{equation}\label{Ydef}
Y =  \frac{\im}{2}( \, |p|^2 \, Z_0 + h\, \bar{p} \, Z_+
   + \bar{h} \, p \, Z_-) + |h|^2\, Z_3,
\end{equation}
so that $\, \langle  Y, \, \bar{Y} \rangle =0$.

$\bullet$ Assume first $\, Y \ne 0$.
One computes from \eqref{struct9}, \eqref{struct90}, that
\begin{equation}
dY + 3\, \phi \, Y = Y_+ \, \frac{1}{2}\, \bomega
+ Y_- \, \frac{1}{2}\, \omega, \n
\end{equation}
where
\begin{align}\label{Ypdef}
Y_+ &=\im(\, p\bar{y}-q\bar{p})\, Z_0 + \im(\, h\bar{y}-\frac{1}{2}p\bar{p})\, Z_+
+\im(\, p\bar{z}-q\bar{h})\, Z_- +  (\,2h\bar{z}-p\bar{h})\, Z_3   \\
Y_- &=\im(\, y\bar{p}-p\bar{q})\, Z_0 + \im(\, z\bar{p}-h\bar{q})\, Z_++
\im(\, y\bar{h}-\frac{1}{2}p\bar{p})\, Z_- +  (\,2z\bar{h}-h\bar{p})\, Z_3,\n
\end{align}
and $\, \hat{X} = [Y]: M \hook \Sigma$ is well defined.
A computation shows
\begin{align}
\langle Y, \, \bar{Y}_{\pm} \rangle &= 0, \, \langle Y_+, \,  \bar{Y}_- \rangle = 0, \n \\
\langle Y_{\pm}, \,  \bar{Y}_{\pm} \rangle &=  2|pz-hy|^2, \n
\end{align}
and, in particular, $\, \hat{X}$ is weakly conformal Legendrian.
Moreover, it is easy to check
$\, \hat{X}$ is in fact Willmore of type \tb{B}
whenever it is an immersion by direct computation.

It thus suffices to show that the zero locus of $\, pz-hy$ is a closed set
with no interior.
Differentiating $\, hq-\frac{1}{2}p^2=0$, we get
$\, hz-py+qz=0$, which implies that
$\, (pz-hy)(qy-px)- \frac{1}{2}(hx-qz)^2=0$.
As in the proof of \tb{Proposition \ref{umbilic}},
\begin{equation}
\partial_{\bar{\xi}}
\begin{pmatrix}
pz-hy \\ hx-qz \\ qy-px
\end{pmatrix} =
\begin{pmatrix}
0 & 1 & 0 \\
\nu(\partial_{\bar{\xi}}) & 0 & 1 \\
0 & \nu(\partial_{\bar{\xi}}) & 0
\end{pmatrix}
\begin{pmatrix}
pz-hy \\ hx-qz \\ qy-px
\end{pmatrix}.\n
\end{equation}
By the same argument as before, either the zeros of $\, pz-hy$ is a closed set
with no interior,
or it vanishes identically. But the latter case
forces $\, \Psi \equiv 0\, $, contrary to our assumption that the surface
is of type \tb{B}.

$\bullet$ Assume $\, Y(m) = 0$ for $\, m \in M$, i.e., $\, h(m) = p(m) = 0$.
Following the notation of \tb{Proposition \ref{umbilic}},
let $(h, \, p, \, q) = \xi^k\, (h_1, \, p_1, \, q_1)$ for a local coordinate
$\, \xi$ centered at $\, m$. There are two cases.

Case $\, (h_1(m), \, p_1(m)) \ne (0, 0)$($k \geq 1$):
Let $\, \lambda = (\xi \, \bar{\xi})^{-k}$, and extend $\, \hat{X}$ across
$\, m$ smoothly by $\, \hat{X} = [ Y_0] = [ \lambda \, Y]$ by using
a local cut off function. We have
\begin{equation}
dY_0 \equiv \lambda \, (Y_+ \, \frac{1}{2}\, \bomega
+ Y_- \, \frac{1}{2}\, \omega), \quad \mod Y_0 \n
\end{equation}
A computation shows
\begin{equation}
pz - hy \equiv k \, \xi^{2k-1}(p_1h_1 - h_1p_1)
\equiv 0 \quad \mod \xi^k. \n
\end{equation}
Since $\, k \geq 1$ by assumption,
$\, \hat{X}$ is smooth across $\, m$.
A similar argument as above
shows that the singular locus of  $\, \hat{X}$
is a closed set with no interior.

Case $\, (h_1(m), \, p_1(m)) =(0, 0)$:
As before, from the choice of $\, k$, $\, q_1(m) \ne 0$,
and $\, h_1 = \frac{1}{2q_1} \, p_1^2$.
Set $\, \lambda = (p \, \bar{p})^{-1}$, and define
$\, \hat{X} = [ Y_0] = [ \lambda \, Y]$ by using a local cut off function.
From \eqref{Ydef}, this gives a well defined extension of $\, \hat{X}$ across
$\, m$. Moreover, since $\, h = \frac{1}{2q} \, p^2$,
it is easy to check that $\, pz - hy \equiv 0 \mod p^2$.
Thus $\, \hat{X}$ is smooth across $\, m$. A similar argument as above
shows that the singular locus of  $\, \hat{X}$ is a closed set with no interior.
\sq






\vsp{2pc}
\section{Willmore Legendrian surfaces with constant dual}

Willmore Legendrian surfaces of type \tb{C} can be completely
characterized in terms of the following set of holomorphic data.
\begin{equation}\label{data}
\left\{
\begin{array}{ll}
\mbox{A compact Riemann surface $\, M$ with a divisor $\, D$.}   \\
\mbox{A meromorphic immersion $f=(f^1, \, f^2): \, M \hook \mathbb{C}^2$
with simple poles on  $\,D$}.   \\
\mbox{Re$\, (f^1 df^2 - f^2 df^1)$ is exact on  $\, M - D$
with zero logarithmic term on $D$.}
\end{array} \right.
\end{equation}
For the definition of \emph{zero logarithmic term},
see the proof of \tb{Theorem \ref{typeC}} below.

Suppose a Willmore Legendrian surface $\, M$  of type \tb{C} is not totally
geodesic. Form the structure equations \eqref{struct90}, we may translate
$\, p =0$ on the dense open subset $\, M^*$ where $\, h \ne 0$. Then
$\, q=0, \, x=0, \, y=0$, and from \eqref{Ydef}, \eqref{Ypdef},
the associated Willmore dual is constant on $\, M^*$,
and hence constant on $\, M$.
\begin{theorem}\label{typeC}
Let $\, X: \,M \hook \Sigma$ be a compact, connected Willmore Legendrian
surface of type \tb{C}. Assume $\, X$ is not totally geodesic, and denote
the constant Willmore dual $\, \hat{X}(M) = p_{\infty}$.
Let $\, D = X^{-1}(p_{\infty})$, which is a divisor on $\, M$.
Let $\, P: \Sigma - \{p_{\infty}\} \to \C^2$
be a stereographic projection from $\, p_{\infty}$.
Then $\, P \circ X: \, M-D \hook \C^2$
is a complete, exact, minimal Lagrangian surface
of finite total curvature with embedded ends of zero logarithmic growth.
Hence, upon a linear change of coordinates of $\, \C^2$,
$\, P \circ X: \, M-D \hook \C^2$ is holomorphic,
and completes across $\, D$ as a holomorphic immersion
$\, f_X: \, M  \hook \C P^2 = \C^2 \cup \C P^1_{\infty}$ with transversal
intersection with $\, \C P_{\infty}^1$ on $\, D$.

Conversely, suppose a Weierstra\ss\, type data \eqref{data} is given.
Then there is an associated conformal Legendrian lift
$\, X_f: M-D \hook \Sigma$ that completes across $\, D$ as
a Willmore Legendrian immersion of type \tb{C}.

In this case, the value of the Willmore functional
$\, \W(X) = 4\, \pi (|D|+g-1)$, where $\, g$ is the genus of $\, M$.
\end{theorem}
\Proof From the geometric description of Willmore dual,
$\, P \circ X: \, M-D \hook \C^2$ is obviously an exact, minimal Lagrangian
surface.
But minimal Lagrangian surfaces in $\, \C^2$ are holomorphic
curves under a linear change of coordinates of $\, \C^2$.
Let $\, f: M-D \hook \C^2$ denote this holomorphic curve.
It is of finite total curvature for
\begin{equation}\label{ChO}
\int_{M-D} \,- K \,dA =  \int_{M-D} \, \frac{1}{2} \, |B_0|^2 dA
      = \int_{M-D} \Omega_X = \W(X)< \infty,
\end{equation}
where $\, K$ is the Gau\ss\; curvature of $\, f(M-D)$.
It is complete with simple poles and embedded ends
because, from the construction of
pseudoconformal stereographic projection in \tb{Section 1},
$\, \E = \langle w_0, \, w_3 \rangle^{\perp}$
can be identified with the contact hyperplane at infinity $\, p_{\infty}$.
That each end is of zero logarithmic growth
follows from the argument for the converse below.

By the classical theorem of Chern-Osserman [ChO], the Gau\ss \;map of
$\, f_X: M-D \hook \C^2$ extends holomorphically across $\, D$.
Since $\, f_X$ is itself holomorphic,
this implies the holomorphic extension
$\, f_X: M  \hook \C P^2 = \C^2 \cup \C P^1_{\infty}$ is well defined.
It is easy to see that $\, f_X$ is an immersion on $\,D$, for the original
$\, X: \,M \hook \Sigma$ is an immersion.

Conversely, suppose such a holomorphic curve
$\, f: M-D \hook \C^2$ is given.
Upon an appropriate linear change of coordinates of $\, \C^2$,
Re$ \, d(f^1 \, df^2 - f^2 \, df^1)$
becomes the standard symplectic form of $\, \C^2$,
and $\, f$ is exact, minimal, and Lagrangian.
Let
\begin{equation}\label{Llift}
X_f: M -D \hook \Sigma = (\R \oplus \C^2) \cup \{ p_{\infty} \}
\end{equation}
be an associated Legendrian lift.
It suffices to show that
$\, X_f$ extends across $\, D$ as a smooth Legendrian
immersion(then it will be necessarily of type \tb{C}).

Take $\, m\in D$, and let $\, \xi$ be a local coordinate of $\, M$ centered
at $\,m$. From the given description, $\, f$ has a local expansion
\begin{equation}
f = f_{-1} \, \frac{1}{\xi} + f_0 + f_1 \, \xi + \, ... \, \n
\end{equation}
Since Re$\,(f^1 \, df^2 - f^2 \, df^1)$ is exact without logarithmic terms,
it cannot have any Re$\, \frac{1}{\xi} \, d\xi$ term
(this is the definition of zero logarithmic terms),
and hence $\,(f^1 \, df^2 - f^2 \, df^1)$ itself
cannot have any $\, \frac{1}{\xi} \, d\xi$ term.
Note that the associated Legendrian lift is not smooth
if there is a nonzero Re $\frac{1}{\xi} \, d\xi$ term.
Thus for meromorphic curves in $\, \C^2$ with simple poles on $\, D$,
each ends has zero logarithmic growth when the meromorphic 1-form
$\, (f^1 \, df^2 - f^2 \, df^1)$ is locally exact.
We may thus write
\begin{equation}
dt = \mbox{Re} \,\big(( \, c_{-1} \xi^{-2} + c_0 + \, ... \,)\, d\xi \big). \n
\end{equation}
for a local smooth function $\, t$ in a neighborhood of $\, m$.
It now follows easily from the inversion formula \eqref{Qform} that
the Legendrian lift $\, X_f$ extends smoothly across $\, m$.

A complete minimal surface of finite total curvature in a Euclidean space
with $\, d$ embedded ends is conformally equivalent to
a compact Riemann surface with
$\, d$ points removed, and with the total Gau\ss\, curvature
$\, -4\pi(d+g-1)$ [ChO]. $\, \W(X) = 4\, \pi (|D|+g-1)$
follows from \eqref{ChO}.
\sq

\Remark $\, D \ne \emptyset$.
\begin{corollary}
Let $\, X: \,S^2 \hook \Sigma$ be a Willmore Legendrian immersion.
Then it is of type \tb{C}, and
$\, \W(X) = 4\, \pi (d-1)$ for some positive integer $\, d$.
\end{corollary}
\Proof By Riemann-Roch, $\, S^2$ does not support any nonzero
holomorphic differential. Hence both $\, \Phi$ and $\, \Psi$ are zero.
\sq

The definition of Legendrian surface can be suitably modified
to include nonorientable surfaces. In this regard, we have the following
analogue of the well known theorem that $\, \R P^2$ cannot occur as
a minimal surface in $\, S^3$, nor as a minimal Legendrian surface in $\, S^5$.
The difference is that in our case it follows from
the global geometry of Willmore Legendrian spheres,
whereas in minimal surface case it follows from the rigidity
of such minimal spheres.
\begin{corollary}
$\, \R P^2$ cannot occur as a Willmore Legendrian surface in $\, \Sigma$.
\end{corollary}
\Proof Let $\, X_0: S^2 \hook \Sigma$ be a Willmore Legendrian immersion
invariant under antipodal involution $\, \tau$ of $\, S^2$,
so that it factors through $\, X=X_0/\tau: \, \R P^2 = S^2 / \tau \hook \Sigma$.
By \tb{Theorem \ref{typeC}}, there exists an associated  holomorphic
immersion $\, \R P^2 \hook \C P^2$, a contradiction.
\sq

\begin{example}\textnormal{
Consider a meromorphic map
$\, f= (\, f^1, \, f^2): \C P^1 \to \C^2$ defined by
\begin{equation}
\begin{pmatrix}
f^1(z) \\ f^2(z)
\end{pmatrix}=
\begin{pmatrix}
\frac{P(z)}{Q(z)} \\
\frac{A(z)}{B(z)}
\end{pmatrix}
\end{equation}
where $\, P, \, Q$ are polynomials of degree $\, k$, and
$\, A, \, B$ are polynomials of degree $\, l$.
By taking generic polynomials, we may assume
$\, f$ is a meromorphic immersion
with $\, k+l$ simple poles on $\, D$.
From
\begin{equation}
f^1df^2-f^2df^1=
\frac{PQ(BA'-AB')-AB(QP'-PQ')}{Q^2B^2} dz, \n
\end{equation}
further generic choice of $\, P, \,Q, \,A, \,B$
implies that
the meromorphic 1-form $\, f^1df^2-f^2df^1$ has no residue,
and is exact on $\, \C P^1 - D$.
By \tb{Theorem \ref{typeC}},
we  have a Willmore Legendrian sphere
$\, X_f: \C P^1 \to \Sigma$ with
$\, \W(X_f) = 4 \pi (k+l-1)\,$.
This is in contrast with the case of
Willmore spheres in conformal 3-sphere,
where certain values $\, 4\pi k\,$ are prohibited [Br3].
It is not difficult to construct $\, f$ so that
$\, X_f$ is an embedding on $\, \C P^1 - D$, [LY].}
\end{example}

A question naturally arises as to how many compact Riemann
surfaces support the Weierstra\ss\, type data \eqref{data}.
\begin{theorem}\label{theorem42}
Every compact Riemann surface $\, M$ admits a generally one to one,
conformal Legendrian immersion $\,X: M \hook \Sigma$
as a Willmore Legendrian surface of type \tb{C}.
\end{theorem}
\Proof
We follow the construction of Bryant [Br1], and the notations therein.
Let $\, (z^0, \, z^1, \, z^2)$ be a coordinate of $\, \C^3$.
Bryant showed that for any compact Riemann surface $\, M$, there exists
a meromorphic embedding $\, F=(F^0, \, F^1, \, F^2): M \hook \C^3$
with simple poles on a divisor $\, D \subset M$ as an integral curve of
the holomorphic contact form $\, dz^0 -z^1dz^2+z^2dz^1$.
Let $\, f=(F^1, \, F^2):M \hook \C^2$, which is necessarily
an exact meromorphic immersion with simple poles on $\, D$.
The construction in [Br1, Theorem G] implies $\,f$ is generally one to one.
The Legendrian lift of $\, f$, \eqref{Llift}, is a graph over $\, f$,
and hence a generally one to one immersion on $\, M - D$.
\sq

The argument in the proof above shows that for every holomorphic
Legendrian immersion $\, F:M \hook \C P^3$, there exist associated
(branched) Willmore Legendrian surfaces of type \tb{C}.
For a discussion of the moduli of holomorphic Legendrian curves, [ChMo].

\begin{corollary}
Every compact Riemann surface can be conformally immersed
in $\, \C^2$ as an exact, algebraic Lagrangian surface.
\end{corollary}

It is known that every compact Riemann surface can be
conformally embedded in Euclidean 3-space as an algebraic surface [Ga].



\vsp{4pc} \noindent
\tb{\large{Appendix.  Euler-Lagrange equation}}

\noindent \emph{Proof of  Proposition \ref{EL}}. $\;$
Let $\, X_t:(-\epsilon, \epsilon) \times M \hook \Sigma$ be one parameter
family of Legendrian immersions with the associated 3-adapted bundle
$\, F_t \to (-\epsilon, \epsilon) \times M$,
and the induced Maurer-Cartan form $\, \pi$ such that
\begin{equation}
\pi = (\pi^A_B) =
\begin{pmatrix}
\phi+(x_0+ \im y_0)dt & \nu_j+\im \mu_j + ( v_j + \im u_j)dt & q_i \omega^i+wdt \\
\omega^i + \im \, \lambda_i \, dt& \pi^i_j & \mu_j+\im \nu_j + ( u_j + \im v_j)dt \\
\lambda \, dt &-\im ( \omega^j - \im \, \lambda_j \, dt ) & -\phi+(-x_0+ \im y_0)dt
\end{pmatrix},\n
\end{equation}
where
\begin{equation}
\pi^i_j -\delta_{ij} \pi^0_0 =
\alpha^i_j+ \im \beta^i_j+(x^i_j+ \im y^i_j)dt, \quad
\pi^i_j + \bar{\pi}^j_i = 0,\n
\end{equation}
and $\, \phi, \, \omega^i, \, \alpha, \, \beta, \, \mu, \, \nu$ are without
$\, dt$-terms.

Let $\, \Omega =$ Im $(\pi^1_1 - \pi^2_2) \w$ Im $\pi^1_2$, and set
\begin{equation}
\Omega_t = \Omega - dt \w(\frac{\partial}{\partial t} \lrcorner \, \Omega). \n
\end{equation}
Then $\, \frac{\partial}{\partial t}  \lrcorner \, \Omega_t = 0$, and
\begin{equation}
\mathcal{L}_{\frac{\partial}{\partial t}} \, \Omega_t |_{t=0}
\equiv \frac{\partial}{\partial t} \lrcorner \, d\Omega|_{t=0}
\quad \mod \mbox{\; exact form}.\n
\end{equation}
It thus suffices to compute
$\, \frac{\partial}{\partial t} \lrcorner \, d\Omega|_{t=0}$.

From $\, -d\pi = \pi \w \pi$, we find
\begin{align}
d\lambda &\equiv \lambda(\pi^0_0 -\pi^3_3)-2\lambda_k \, \omega^k,  \n\\
d\lambda_i &\equiv -\lambda_j \, \omega^i_j+y^i_j \, \omega^j-\lambda \, \nu_i,
 \quad \mod \; dt. \n
\end{align}
A direct computation with these equations gives
\begin{equation}
\frac{\im}{2} \frac{\partial}{\partial t} \lrcorner \, d\Omega|_{t=0}
= \, y^i_j \, p_{ij} \, \omega^1 \w \omega^2 +
(\lambda_1\nu_1-\lambda_2\nu_2)\w \beta^1_2
-\beta^1_1 \w (\lambda_1\nu_2+\lambda_2\nu_1), \n
\end{equation}
where $\, \mu_i = p_{ij} \, \omega^j$.
We wish to put this expression into a form involving $\, \lambda$ only
by integration by parts. After a short computation,
\begin{equation}
\frac{\im}{2} \frac{\partial}{\partial t} \lrcorner \, d\Omega|_{t=0}
+ d(\lambda_1 \, \mu_2 - \lambda_2 \, \mu_1)
= (\lambda_1q_1+\lambda_2q_2)-\lambda(\nu_1\w \mu_2-\nu_2\w \mu_1). \n
\end{equation}
Adding $\, \frac{1}{2}\, d(\lambda(q_1\omega^2-q_2\omega^1))$, we get
\begin{equation}
\frac{\im}{2} \frac{\partial}{\partial t} \lrcorner \, d\Omega|_{t=0}
+ d(\lambda_1 \, \mu_2 - \lambda_2 \, \mu_1)
+ \frac{1}{2}\, d(\lambda(q_1\omega^2-q_2\omega^1))
= \lambda \, r \, \omega^1\w\omega^2. \n
\end{equation}
Since $\, \lambda$ is an arbitrary function of compact support,
this forces $\, r=0$.
\sq

\vspace{4pc}

\noindent
\begin{center}
\textbf{\large{References}}
\end{center}

\noindent [BB]Babich, M.; Bobenko, A., \emph{
Willmore tori with umbilic lines and minimal surfaces in hyperbolic space},
Duke Math. J.  72  (1993),  no. 1, 151--185

\noindent [Br1]  Bryant, Robert L., \emph{
Conformal and minimal immersions of compact surfaces into the $4$-sphere},
J. Differential Geom.  17  (1982), no. 3, 455--473

\noindent [Br2]  $\underline{\qquad}$, \emph{
A duality theorem for Willmore surfaces},
J. Differential Geom.  20  (1984),  no. 1, 23--53

\noindent [Br3]  $\underline{\qquad}$, \emph{
Surfaces in conformal geometry},
The mathematical heritage of Hermann Weyl,
Proc. Sympos. Pure Math. 48 (1988), 227--240

\noindent [ChM] Chern, S. S.; Moser, J. K., \emph{
Real hypersurfaces in complex manifolds},
Acta Math.  133  (1974), 219--271

\noindent [ChO] $\underline{\qquad}$; Osserman, Robert, \emph{
Complete minimal surfaces in euclidean $n$-space},
J. Analyse Math.  19  (1967), 15--34

\noindent [ChMo] Chi, Quo-Shin; Mo, Xiaokang, \emph{
The moduli space of branched superminimal surfaces of a fixed degree,
genus and conformal structure in the four-sphere},
Osaka J. Math. 33 (1996), no. 3, 669--696

\noindent [Ga] Garsia, Adriano M., \emph{
On the conformal types of algebraic surfaces of euclidean space},
Comment. Math. Helv. 37 (1962/1963), 49--60.

\noindent [Gro] Gromov, M., \emph{
Pseudoholomorphic curves in symplectic manifolds},
Invent. Math. 82 (1985), no. 2, 307--347

\noindent [He] Helein, F., \emph{
Willmore immersions and loop groups},
J. Differential Geom.  50  (1998),  no. 2, 331--385

\noindent [LY] Li, Peter; Yau, S. T.,\emph{
A new conformal invariant and its applications to the Willmore conjecture
and the first eigenvalue of compact surfaces},
Invent. Math.  69  (1982), no. 2, 269--291

\noindent [Min] Minicozzi, William P., II, \emph{
The Willmore functional on Lagrangian tori: its relation to area
and existence of smooth minimizers},
J. Amer. Math. Soc. 8 (1995), no. 4, 761--791

\noindent [MU]
Montiel, Sebastian; Urbano, Francisco, \emph{
A Willmore functional for compact surfaces in the complex projective plane},
J. Reine Angew. Math. 546 (2002), 139--154

\noindent [Ro] Ros, Antonio, \emph{
The Willmore conjecture in the real projective space},
Math. Res. Lett.  6  (1999),  no. 5-6, 487--493

\noindent [To] Topping, P., \emph{
Towards the Willmore conjecture},
Calc. Var. PDE, 11 (2000), no. 4, 361--393

\vspace{2pc}

\noindent

\noindent Sung Ho Wang \\
\noindent Department of Mathematics \\
\noindent Kias \\
\noindent Seoul, Corea 130-722 \\
\texttt{shw@kias.re.kr}

\end{document}